\documentclass[times,review,sort&compress]{elsarticle}
\usepackage[active]{srcltx}
\usepackage{bbm}
\usepackage{mathrsfs}
\usepackage{tipa}
\usepackage{amsthm}
\usepackage{amssymb}
\usepackage{stackrel}
\usepackage{amsmath,amscd}
\newtheorem{dl}{Theorem}[section]
\newtheorem{tl}[dl]{Corollary}
\newtheorem{yl}[dl]{Lemma}

\newtheorem{lz}[dl]{Example}

\newtheorem{remark}[dl]{Remark}

\numberwithin{equation}{section}
\newproof{pot333}{Proof of Proposition \ref{ylthird}}
\newcommand{\poq}[2]{(#1;q)_{#2}}

\def\qed{\hfill \rule{4pt}{7pt}}
\def\pf{\noindent {\it Proof.} }
\newcommand{\poqq}[2]{(#1;q^2)_{#2}}

\begin{document}
\title{On Andrews--Warnaar's   identities of  partial theta functions}
\author{Jin Wang\fnref{fn1}}
\fntext[fn1]{E-mail address: jinwang2016@yahoo.com}
\author{Xinrong Ma\fnref{fn2,fn3}}
\fntext[fn2]{Corresponding author.}
\fntext[fn3]{E-mail address: xrma@suda.edu.cn.}
\address[P.R.China]{Department of Mathematics, Soochow University, SuZhou 215006, P.R.China}
\markboth{J. Wang and  X. Ma}{On Andrews--Warnaar's   identities of  partial theta functions}
\begin{abstract}
In this paper we set up a bivariate representation of partial theta functions which not only unifies some famous identities for partial theta functions due to Andrews and Warnaar, et al. but also unveils a new characteristic of such identities. As  further applications, we establish a general form of Warnaar's identity  and a general $q$--series transformation associated with Bailey pairs via the use of the power series expansion of partial theta functions.
\end{abstract}
\maketitle
\vspace{20pt}
\parskip 7pt
\section{Introduction}
Throughout this paper, we adopt the standard notation and terminology for  $q$-series from the book  \cite{10} (Gasper and Rahman, 2004).  As customary, the $q$--shifted factorials of complex variable $x$ with the base $q$ are  given by
\begin{eqnarray*}
(x;q)_\infty
&=&\prod_{n=0}^{\infty}(1-xq^n)\qquad\mbox{and}\\
(x;q)_n&=&\frac{(x;q)_\infty}{(xq^n;q)_\infty}
\end{eqnarray*}
for all integers $n$. For integers $m\geq 1$, we  employ
 the multi-parameter compact notation
\[(a_1,a_2,\ldots,a_m;q)_\infty=(a_1;q)_\infty (a_2;q)_\infty\ldots (a_m;q)_\infty.\]
In \cite{10},  the  $q$--Heine ${}_2\phi_1$ series with the
base $q$ and the argument $x$ is defined by
\begin{eqnarray*}
{}_{2}\phi _{1}\left[\begin{matrix}a,&b\\ &c\end{matrix}
; q, x\right]&=&\sum _{n=0} ^{\infty }\frac{\poq {a,b}{n}}{\poq
{c}{n}}\frac{x^{n}}{\poq{q}{n}}.
\end{eqnarray*}
Sums of the form
\begin{eqnarray}
\sum_{n=0}^\infty
q^{An^2+Bn}x^{n}\label{theta-sum}
\end{eqnarray}
 are called partial theta functions owing to the fact that
 \begin{eqnarray*}
\sum_{n=-\infty}^\infty
q^{An^2+Bn}x^{n}\quad  (x\neq 0)
 \end{eqnarray*}
 are often referred to as  (complete)  theta functions. For the complete theta functions, there holds the famous  Jacobi
triple product identity (\cite[(II.28)]{10})
\begin{eqnarray}
 \sum_{n=-\infty}^\infty (-1)^nq^{n(n-1)/2}x^n=
 \poq{q,x,q/x}{\infty}.\label{triple}
 \end{eqnarray}
Following \cite{andrews, andrews-31},  any $q$--series identity containing sums
such  as (\ref{theta-sum}) is called  a partial theta function identity.   Partial theta function identities first appeared in Ramanujan's   legendary lost notebook \cite{raman-s}, wherein he recorded many identities of such sort  without any proofs. A full treaties together with  a more complete bibliography on this topic can be found in
\cite[Chap. 6]{andrews1} by G. E. Andrews and B. C. Berndt. Here, we choose from \cite{andrews1} some typical identities as  illustrations.
\begin{yl}
\begin{description}
            \item[(i)] \mbox{\rm (cf. \cite[p.12]{raman-s} or \cite[Entry 6.6.1]{andrews1})}
\begin{align}
  \sum_{n=0}^\infty q^{n(n+1)}a^n=\sum_{n=0}^\infty \frac{\poq{q^{n+1}}{n}q^n}{\poq{a}{n+1}\poq{q/a}{n}}-\sum_{n=0}^\infty
 \frac{a^{3n+1}q^{n(3n+2)}(1-aq^{2n+1})}{\poq{a,q/a}{\infty}}. \label{war1}
 \end{align}
\item[(ii)]\mbox{\rm (cf. \cite[p.4; p. 29]{raman-s} or \cite[Entries  6.3.9 and 6.3.11]{andrews1})}
\begin{align}
 \sum_{n=0}^\infty q^{n(n+1)/2}a^n&=\sum_{n=0}^\infty \frac{\poqq{q}{n}q^n}{\poq{a}{n+1}\poq{q/a}{n}}+\sum_{n=0}^\infty
 \frac{(-1)^{n+1}a^{2n+1}q^{n(n+1)}}{\poq{-q,a,q/a}{\infty}} \label{war2}\\
 &=\sum_{n=0}^\infty \frac{\poqq{q}{n}q^{2n}}{\poqq{a}{n+1}\poqq{q^2/a}{n}}+\sum_{n=0}^\infty
 \frac{(-1)^{n+1}a^{3n+1}q^{n(3n+2)}(1+aq^{2n+1})}{\poq{-q}{\infty}\poqq{a,q^2/a}{\infty}}. \label{war3}
\end{align}
\end{description}
\end{yl}
According to all publications accessible to us, we see that the first systematic study on partial theta function identities is due to  G. E. Andrews  \cite{andrews,andrews-31}. Since then, particularly in a very recent decade, much more progress on partial theta functions  and allied identities has been made.  All that, in turn, show that the partial theta function (\ref{theta-sum})  is really important and worth of study.

By a short excursion on the methods and results published in the research literature, we feel that the study of partial theta functions mainly covers  over three  kinds of research areas: analysis, combinatorics, and $q$--series.

 As is to be expected, analytic area is devoted to, when (\ref{theta-sum}) is viewed as an analytic function of variable $x$ with $q$ being a parameter, asymptotic expansions, zeros, expansion coefficients etc. For asymptotic expansions, we refer the reader  to  B.C. Berndt and B. Kim \cite{new-1},  S. Jo and B. Kim \cite{new-8}, B. Kim and J. Lovejoy \cite{new-9}, R. R. Mao \cite{new-6},  R. J. McIntosh \cite{new-11}, B. Kim,  E. Kim, and J. Seo \cite{new-92}. Especially noteworthy is that  in \cite{new-8}, Jo and  Kim
  defined the integer sequence $\{a_n(d)\}_{n\geq 0}$  by the expansion
 \begin{eqnarray*}
 \sum_{n=0}^\infty  (-1)^nq^{(dn-d+2)n/2}=\poq{q}{\infty}\sum_{n=0}^\infty a_n(d)q^n
 \end{eqnarray*}
and discovered many interesting combinatorial and analytic properties of $a_n(d)$.  Besides,
Bringmann and Folsom  confirmed in \cite{new-122}  a conjecture of Berndt and Kim (cf.\cite[Conjecture 1]{new-1}), which claims that  for sufficiently large $n$, the
coefficients $a_n$ in the asymptotic expansion (corresponding to   (\ref{theta-sum}) with $A=-x=1$ and $B$ being  positive integers)
\begin{eqnarray*}
2\sum_{n=0}^\infty(-1)^n
\bigg(\frac{1-t}{1+t}\bigg)^{n^2+Bn}=\sum^\infty_{n=0}a_nt^n
\end{eqnarray*}
 have the same sign.

For recent works on zeros of a partial theta function,  the interested
reader is referred to the earlier paper \cite{andrews-31} by G. E. Andrews, and  a series of papers   \cite{new-32,new-35,new-34,new-344,new-33,new-31,new-36} by V. P. Kostov,   as well as  the latest papers \cite{new-12} by T. Prellberg and   \cite{new-4} by A. Sokal. By considering  the partial theta function (\ref{theta-sum})  as an entire function in $x$, Andrews and Kostov independently  investigated  the double and multiple zeros, as well as the spectrum of the partial theta functions.  Among the various results, the most noteworthy  ones  are Kostov--Shapiro's resolution of the Hardy--Petrovitch--Hutchinson
problem in \cite{new-33} and Kostov's conclusions in \cite{new-36, new-344} that for any fixed $q$,   the partial theta function
 has at most finitely many multiple zeros; and if $y(q)$ is a multiple zero of the function (\ref{theta-sum}) with $x$ replaced by $-xq$, then $|y(q)|<8^{11}$.
 It should be mentioned that  in \cite{new-4}, Sokal showed that if $x_0(q)$ denotes the leading root of partial theta function (\ref{theta-sum}) with $A=B=1$, then the coefficients   in the expansions of  $x_0(q), -1/x_0(q), 1/x_0^2(q)$ as power series in $q$ are strictly negative except for few initial terms.

 As far as combinatorial area is concerned, there is an amount of  literature on combinatorial bijective  proofs and number--theoretic interpretations  of partial theta function identities. Some recent literature which we  would like to recommend  include  G. E.  Andrews \cite{andrews21},   K. Alladi and A. Berkovich \cite{berk0}, K. Alladi \cite{all-new-11, all-new-12}, B. C. Berndt, B. Kim, and A. J. Yee \cite{new-7},  B. Kim and J. Lovejoy  \cite{new-9}, B. Kim  \cite {new-91}, K. Q. Ji, B. Kim, and J. S. Kim  \cite{ji}, A. J. Yee \cite{new-10}, T. Prellberg \cite{new-12}. Perhaps most noteworthy in this respect is that Andrews pointed out in \cite{andrews21} via the ninth and tenth open problems that the partial theta function identity from Ramanujan's lost notebook \cite[Enrty 1.6.2]{andrews1}
 \begin{align}
  \sum_{n=0}^\infty q^{n^2}a^n=1+\sum_{n=1}^\infty \frac{\poq{-q}{n-1}}{\poqq{-aq^2}{n}}a^nq^{n(n+1)/2} \label{andrews-yee}
 \end{align}
 has a connection with some parity problems in integer partitions. Shortly after \cite{andrews21}, these two problems were solved by Yee in \cite[Theorem 2.2]{new-10} along with a combinatorial proof of (\ref{andrews-yee}). We remark in passing that, as \cite{new-9,ji} display, the expansion  coefficients  given by some partial theta function identities are  unimodal.

 The third area of research is  concerned about possible relationship of partial theta function identity with Bailey pairs and $q$--series.  The reader may  consult B. Kim and J. Lovejoy \cite{new-13,new-93}, K. Bringmann, A. Folsom, and R. C. Rhoades \cite{new-5},  J.  Lovejoy \cite{new-2},  A. Schilling and S.  O.  Warnaar \cite{warnaar00}, S.  O.  Warnaar \cite{warnaar} for further information. We especially refer the reader to \cite{new-5} for the connection between partial theta functions and mock theta functions. Here, what we want to emphasize  is that   in order to gain deeper insights into ``hidden" structure of (\ref{war1})--(\ref{war3}),  Warnaar set up  in his remarkable paper \cite{warnaar} the following partial theta function identity.
 \begin{yl}[Warnaar \mbox{\cite[Theorem 1.5]{warnaar}}]\label{yl2} For any complex number $x$, define the partial theta function
\begin{eqnarray}
\theta(q,x)=\sum_{n=0}^\infty \tau(n)x^n,\label{thetafun}
\end{eqnarray}
where $\tau(n)=(-1)^nq^{n(n-1)/2}.$ Then for any  $a,b$, it holds
  \begin{eqnarray}
\theta(q,a)+\theta(q,b)-1=\poq{q,a,b}{\infty}\sum_{n= 0}^\infty \frac{\poq{ab/q}{2n}}{\poq{q,a,b,ab}{n}}q^n.\label{2}
\end{eqnarray}
\end{yl}
Evidently, (\ref{2}) is a generalization of Jacobi's triple product identity  (\ref{triple}).  More importantly, by virtue of (\ref{2}), Warnaar established a $q$--series  transformation formula which can be restated as  follows.
   \begin{yl}\label{warnaar-bailey}\mbox{\rm (Warnaar \cite[Corollary 4.1]{warnaar})}
\begin{eqnarray}
  &&\sum_{n=0}^\infty \frac{\poq{q}{2n}q^n}{\poq{a}{n+1}\poq{q/a}{n}}\beta_n(q)=
  (1-q)\sum_{n=0}^\infty \frac{(-1)^nq^{-n(n-1)/2}a^n}{1-q^{2n+1}}\alpha_n(q)\label{war4}\\
  &&\qquad\qquad+\frac{1}{\poq{q^2,a,q/a}{\infty}}\sum_{k=1}^\infty (-1)^{k+1}a^kq^{k(k-1)/2}
 \sum_{n=0}^\infty\frac{q^{(1-k)n}(1-q^{k(2n+1)})}{1-q^{2n+1}}\alpha_n(q)\nonumber
 \end{eqnarray}
provided all sums converge. Hereafter,  the Bailey pair  $(\alpha_n(t),\beta_n(t))$ relative to the parameter $t$ is given by
\begin{eqnarray}
\beta_n(t)=\sum^{n}_{k=0} \frac{\alpha_k(t)}
{\poq{q}{n-k}\poq{tq}{n+k}}.\label{defnew}
 \end{eqnarray}
\end{yl}
As it turns out, (\ref{war4}) serves as a bridge connecting partial theta function (residual) identities with Bailey pairs. An  excellent survey on Bailey pairs and their applications via Bailey's lemma and Bailey chains can be   found in \cite{warnaar0}.

 Suitable choices of $(\alpha_n(q),\beta_n(q))$ in (\ref{war4}) led Warnaar  to
 simple proofs of identities (\ref{war1})--(\ref{war3}).  For the details the reader may see \cite[Section 3.1]{warnaar}. As further developments of \cite{warnaar}, there are  some recent works deserving our special mention. At first, continuing in the line of \cite{warnaar}, Lovejoy \cite{new-2} explored extensively  how to derive  conjugate Bailey pairs from residual identities of Ramanujan--type partial theta identities. Second,  in \cite{new-93,new-13} Kim and Lovejoy  extended Warnaar's results in \cite{warnaar} to those in terms of partial and ternary indefinite  theta functions. Third,  Ji and Zhao  \cite{ji-new} used the Bailey transform and the conjugate Bailey pair of Warnaar  to show Garvan's two--variable Hecke--Rogers identities for the universal mock theta
functions (also see \cite{garvan,new-5}).  These works  convince us  that  Kim--Lovejoy's and Warnaar's methods pave a new way to Bailey pairs and Bailey chains.

   Let us back to Warnaar's identity (\ref{2}).  In their another paper  \cite{warnaar1}, Andrew and Warnaar  discovered  a product version of  (\ref{2}), namely (\ref{1}) as below. For our purpose, we now restate it by the following lemma.

\begin{yl}[Andrews--Warnaar \mbox{\cite[Theorem 1.1]{warnaar1}}]\label{yl2-12} With the same notation as above.  Then for any complex numbers $a,b$, it holds
  \begin{eqnarray}
\theta(q,a)\theta(q,b)=\poq{q,a,b}{\infty}\sum_{n=0}^\infty \frac{\poq{ab/q}{2n}}{\poq{q,a,b,ab/q}{n}}q^n.\label{1}\end{eqnarray}
\end{yl}

As further discussion on  identities (\ref{2}) and (\ref{1}),  Berkovich pointed out in \cite{berv}  that (\ref{1}) is a special case of the Gasper--Rahman product formula \cite[p.235,(8.8.18)]{10}. He also showed that (\ref{2}) is equivalent to the following

\begin{yl}[Schilling--Warnaar \mbox{\cite[Lemma 4.3]{warnaar00}}]\label{yl3} Let $\theta(q,x)$ be  given by (\ref{thetafun}). Then for any complex numbers $a,b$, it holds
 \begin{align}
  &&\frac{\theta(q,a)-\theta(q,b)}{a-b}=-\poq{q,aq,bq}{\infty}\sum_{n= 0}^\infty \frac{\poq{ab}{2n}}{\poq{q,aq,bq,ab}{n}}q^n.\label{b4}\end{align}
    \end{yl}
It is also noteworthy that in \cite{berk0}, Alladi and Berkovich established another two--parameter generalization of Jacobi's triple product identity.

\begin{yl}[Alladi--Berkovich \mbox{\cite[Eq.(5.1)]{berk0}}] For the partial theta function $\theta(q,x)$, it holds
  \begin{eqnarray}
  &&\sum_{n=0}^\infty\tau(n)q^{2n}\left\{\frac{a}{a-1}\theta(q,aq^{1+n})+\frac{b}{b-1}\theta(q,bq^{1+n})\right\}
  \frac{\poq{ab}{n}}{\poq{q}{n}}\label{44}
  \\&+&\frac{1-ab}{(1-a)(1-b)}\sum_{n=0}^\infty\tau(n)q^{2n}\theta(q,q^{1+n})\frac{\poq{ab}{n}}{\poq{q}{n}}\nonumber
  =(aq,bq,q;q)_{\infty}.
     \end{eqnarray}
\end{yl}

For more details about (\ref{2}) and (\ref{1}), the reader might consult \cite{andrews,andrews1,warnaar1,berv,10,warnaar} for alternative proofs of using Heine's first and second transformations. In \cite{ma}, the author put forward the so called $t$--coefficient method to deal with all these identities. As for  (\ref{44}),  the reader may consult \cite{berk0} for a pure combinatorial proof in style of integer partitions. In this paper, as one of main purposes, we shall establish a general identity which may serve as common source of all aforementioned theta function identities.
 \begin{dl}\label{maindl}For any integers $m,n\geq 0$ and  complex numbers $a,b$, define
   \begin{eqnarray}
U_m(b)&=&\sum_{k=0}^{\infty}\frac{\poq{q^{1-m}}{k}}
{\poq{q}{k}\poq{bq^{k}}{m}}(1-bq^{2k})b^{2k}q^{2k^2-k+mk},\label{identitiy-1}
 \\
V_{m,n}(a,b)&=&{}_{2}\phi _{1}\left[\begin{matrix}q^{-m},&q^{-n}\\ &abq^{n-1}\end{matrix}
; q, bq^{m+n}\right].\label{identitiy-2}
  \end{eqnarray}
  Then
  \begin{eqnarray}
U_m(b)\theta(q,a)
 =\poq{q,a,b}{\infty}\sum_{n=0}^\infty\frac{\poq{abq^{n-1}}{n}q^n}{\poq{q,a}{n}}
 \frac{V_{m,n}(a,b)}{\poq{b}{m+n}}.\label{mainid}
 \end{eqnarray}
  \end{dl}

Of all instances covered by (\ref{mainid}), the case $m=1$ is the most important one because it sheds light on certain hidden essence of all aforementioned identities of partial theta functions.  In this sense, it deserves a separate theorem.
\begin{dl}[Bivariate representation of partial theta functions]\label{coro000} For any complex numbers $a,b$, define
    \begin{eqnarray}
 L(a,b)=\poq{q,a,b}{\infty}
 \sum_{n=0}^\infty\frac{\poq{ab/q^2}{2n}}{\poq{q,a,b,ab/q^2}{n}}q^{n}.\label{important}
   \end{eqnarray}
 Then we have
 \begin{eqnarray}
 \theta(q,a)=L(a,b)+bL(aq,bq).\label{000}
   \end{eqnarray}
\end{dl}

Curious as it seems, (\ref{000}) can be regarded as a bivariate form of the partial theta function $\theta(q,a)$ despite $\theta(q,a)$ itself is not. It is this expression that leads us to a unified proof for  identities (\ref{2}) and (\ref{b4}). Meanwhile, in a similar argument to that used by Warnaar in proving Lemma \ref{warnaar-bailey},  we are able to  deduce from (\ref{important}) that
\begin{dl}\label{coro111} Let $(\alpha_n(t),\beta_n(t))$ and $L(a,b)$ be given by (\ref{defnew}) and   (\ref{important}), respectively.   Then we have
 \begin{eqnarray}
 \sum_{n=0}^\infty L(aq^{n+1},bq^{n+1})q^n\alpha_n\big(ab/q\big)   =\poq{q,aq,bq}{\infty}\sum_{n=0}^\infty
  \frac{\poq{ab}{2n}}
 {\poq{aq,bq}{n}}q^{n}\beta_n\big(ab/q\big).\label{111-added}
   \end{eqnarray}
\end{dl}

Furthermore, in  the same argument as  Theorem \ref{coro000}  we can derive from the special case $m=2$ of Theorem \ref{maindl} that
\begin{dl}\label{maindl-tl2} For any complex numbers $a,b$, define
    \begin{eqnarray}
 P(a,b)=\poq{q,a,b}{\infty}\sum_{n=0}^\infty\frac{\poq{ab/q^3}{2n}}{\poq{q,a,b,ab/q^3}{n}}q^{n}.\label{important-m2}
   \end{eqnarray}
Then we have
\begin{eqnarray}
\theta(q,a)
 =\frac{q}{q+b}P(a,b)+\frac{b(1+q)}{q+b}P(aq,bq)+\frac{b^2q}{q+b}P(aq^2,bq^2).\label{mainid-new}
 \end{eqnarray}
  \end{dl}

The remainder of our paper is organized as follows.  In  the subsequent section, we shall give the proofs of these theorems (except Theorem \ref{coro111}).  As  direct applications of Theorems \ref{maindl} and \ref{coro000}, Section 3 is devoted to a unified proof of Lemmas \ref{yl2}, \ref{yl2-12}, and \ref{yl3} while Section 4 is about further applications of $L(a,b)$. That is, we shall start with   the power series expansion of $L(a,b)$ to establish a general version of Warnaar's identity (\ref{2}) and  prove Theorem \ref{coro111}.

 \section{The proofs of main results}
We now proceed with Theorems \ref{maindl}, \ref{coro000}, and \ref{maindl-tl2}  whose proofs require the following transformation formula.
 \begin{yl}\label{yl5} For $|tq/(AB)|<1$, it always holds
\begin{eqnarray}
&&\frac{\poq{t,tq/(AB)}{\infty}} {\poq{tq/A,tq/B}{\infty}}\sum_{n=0}^{\infty}\frac{\poq{A,B}{n}}{\poq{q,c}{n}}
\bigg(\frac{tq}{AB}\bigg)^n\label{known}\\
&=&\sum_{n=0}^{\infty}\frac{\poq{t,A,B,tq/c}{n}}
{\poq{q,tq/A,tq/B,c}{n}}(1-tq^{2n})\left(\frac{ct}{AB}\right)^nq^{n^2}.\nonumber
\end{eqnarray}
\end{yl}
\pf It is a special case of Ex. 2.22 in \cite{10} under the specifications
\[(a,c,d,e)\to (t,tq/c,A,B)\]
and then let $b$ tend to infinity.
\qed
\begin{tl}
\begin{eqnarray}
\poq{t}{\infty}\sum_{n=0}^\infty\frac{q^{n^2}t^n}{\poq{q,c}{n}}=\sum_{n=0}^{\infty}\frac{\poq{t,tq/c}{n}}
{\poq{q,c}{n}}(1-tq^{2n})\left(ct\right)^nq^{2n^2-n},&&\label{main}\\
\poq{t}{\infty}\sum_{n=0}^\infty\frac{q^{n^2}t^n}{\poq{q,t}{n}}=\sum_{n=0}^{\infty}\tau(n)t^n.&&\label{main-added}
\end{eqnarray}
 \end{tl}
\pf  Actually  (\ref{main}) is an immediate consequence of taking $A,B\to \infty$ in (\ref{known}) while   (\ref{main-added}) is the special case $c=t$ of (\ref{main}).
\qed
\\

By use of  (\ref{main}) and in a similar argument as  in \cite{warnaar1}, we can set up

\begin{yl}\label{maindl-yl} For any complex numbers $a,b,c$,
  define
   \begin{eqnarray}
 f(b,c)&=&\sum_{n=0}^{\infty}\frac{\poq{b,bq/c}{n}}
{\poq{q,c}{n}}(1-bq^{2n})\left(bc\right)^nq^{2n^2-n},\label{identitiy-1-1}
 \\
 g_n(a,b,c)&=&{}_{2}\phi _{1}\left[\begin{matrix}q^{-n},&q^{1-n}/a\\ &c\end{matrix}
; q, abq^{2n-1}\right].\label{identitiy-2-2}
  \end{eqnarray}
  Then we have
  \begin{eqnarray}
 f(b,c)\theta(q,a)
 =\poq{q,a,b}{\infty}\sum_{n=0}^\infty\frac{q^n}{\poq{q,a}{n}}g_n(a,b,c).\label{identitiy-3-3}
 \end{eqnarray}
 \end{yl}
\pf
   Given $f(b,c)$, according to (\ref{main}), it is very clear that
   \begin{eqnarray}
f(b,c)=\poq{b}{\infty}\sum_{n=0}^\infty
\frac{q^{n^2}b^n}{\poq{q,c}{n}}.\label{oneone}\end{eqnarray}
On the other hand, it is easy to check (or see \cite[(2.1a)]{warnaar1}) that
\begin{eqnarray}
 \sum_{n=0}^\infty\tau(n)a^n=\poq{q,a}{\infty}\sum_{n=0}^\infty
 \frac{q^n}{\poq{q,a}{n}}. \label{twotwo}
 \end{eqnarray}
In the sequel,   by multiplying (\ref{oneone}) with (\ref{twotwo})  we obtain
 \begin{eqnarray*}
 f(b,c)\sum_{n=0}^\infty\tau(n)a^n&=&\poq{q,a,b}{\infty}\sum_{i,j\geq 0}
 \frac{q^j}{\poq{q,a}{j}}\frac{q^{i^2}b^i}{\poq{q,c}{i}}\\
 &=&\poq{q,a,b}{\infty}\sum_{n=0}^\infty  q^n\sum_{i=0}^n
 \frac{\tau^2(i)b^i}{\poq{q,c}{i}}\times \frac{1}{\poq{q,a}{n-i}}.
  \end{eqnarray*}
 Using the basic relation (see \cite[(I.10)]{10})
 \[\poq{x}{n-i}=\frac{\poq{x}{n}}{\poq{q^{1-n}/x}{i}}\frac{\tau(i)q^{i-ni}}{x^i}\]
 to simplify the factors $\poq{q,a}{n-i}$,  we obtain
 \begin{eqnarray*}
 f(b,c)\sum_{n=0}^\infty\tau(n)a^n
 =\poq{q,a,b}{\infty}\sum_{n=0}^\infty\frac{q^n}{\poq{q,a}{n}}g_n(a,b,c),
 \end{eqnarray*}
 where
 \[g_n(a,b,c)=\sum_{i=0}^n
 \frac{\poq{q^{-n},q^{1-n}/a}{i}}{\poq{q,c}{i}}
 (abq^{2n-1})^i.\]
 Note that $g_n(a,b,c)$ just constitutes a $q$--Heine ${}_2\phi_1$  series as given by (\ref{identitiy-2-2}). Therefore, the lemma is proved.
 \qed

Now we are ready to show Theorem \ref{maindl}.

\pf It suffices to take $c=bq^m(m\geq 0)$ in Lemma \ref{maindl-yl}  and assume
\begin{eqnarray*}
f(b,bq^m)=U_m(b)\poq{b}{m}.
\end{eqnarray*}
Then, at first, we have
\begin{eqnarray*}
U_m(b)=\frac{f(b,bq^m)}{\poq{b}{m}}&=&\sum_{n=0}^{\infty}\frac{\poq{b,q^{1-m}}{n}}
{\poq{q}{n}\poq{b}{m+n}}(1-bq^{2n})b^{2n}q^{2n^2-n+mn}\\
&=&\sum_{n=0}^{\infty}\frac{\poq{q^{1-m}}{n}}
{\poq{q}{n}\poq{bq^n}{m}}(1-bq^{2n})b^{2n}q^{2n^2-n+mn}.
\end{eqnarray*}
In the meantime, we have
\begin{eqnarray*}
g_n(a,b,bq^m)&=&{}_{2}\phi _{1}\left[\begin{matrix}q^{-n},&q^{1-n}/a\\ &bq^{m}\end{matrix}
; q, abq^{2n-1}\right]\\
&=&\frac{\poq{abq^{n-1}}{n}}{\poq{bq^{m}}{n}}{}_{2}\phi _{1}\left[\begin{matrix}q^{-m},&q^{-n}\\ &abq^{n-1}\end{matrix}
; q, bq^{m+n}\right].
\end{eqnarray*}The last equality comes from Heine's second transformation of $_2\phi_1$ series [6, (III. 2)].
All these reduces (\ref{identitiy-3-3}) to (\ref{mainid}), which is what we wanted to prove.
\qed

Next are the proofs of Theorems \ref{coro000} and \ref{maindl-tl2}  which seem more or less marvellous.

\pf For $m=1$, we easily deduce from  (\ref{identitiy-1}) that $U_1(b)=1$ and
\begin{eqnarray*}
 V_{1,n}(a,b)=\frac{1-abq^{n-1}+b(1-q^n)}{1-abq^{n-1}}.
  \end{eqnarray*}
   Then
  \begin{eqnarray}
  \theta(q,a)
 =\poq{q,a,b}{\infty}\sum_{n=0}^\infty\frac{\poq{ab/q}{2n}q^n}{\poq{q,a}{n}}
 \frac{1-abq^{n-1}+b(1-q^n)}{\poq{b,ab/q}{1+n}}.\label{22}
 \end{eqnarray}
 Splitting the series on the right--hand side (in short, RHS)  of (\ref{22}) into two parts in such a way that
     \begin{eqnarray*}
 \sum_{n=0}^\infty\tau(n)a^n
 &=&\poq{q,a,bq}{\infty}\sum_{n=0}^\infty\frac{\poq{ab/q}{2n}q^{n}}{\poq{q,a,bq,ab/q}{n}}\\
 &+&bq\poq{q,aq,bq^2}{\infty}\sum_{n=0}^\infty\frac{\poq{abq}{2n}q^{n}}{\poq{q,aq,bq^2,abq}{n}}.
   \end{eqnarray*}
  Taking account of the fact that the partial theta function on the left--hand side (in short, LHS)  of the last  identity is independent of $b$, we can reformulate this identity in symmetric form via the replacement of  $b$ with $bq^{-1}$. The desired identity (\ref{000})  follows.

Analogously, for $m=2$, it is easily found that $U_2(b)=1+bq$ and
\begin{eqnarray*}
V_{2,n}(a,b)=1+b(1+q)\frac{1-q^n}{1-abq^{n-1}}+b^2q^2\frac{(1-q^{n-1})(1-q^n)}{(1-abq^{n-1})(1-abq^{n})}.
\end{eqnarray*}
Substituting these facts into (\ref{mainid}) and then replacing $b$ with $bq^{-2}$, we immediately obtain (\ref{mainid-new}).  Now Theorems \ref{coro000} and \ref{maindl-tl2}  have been proved.
\qed
\section{Unified proofs of partial theta function identities}
As planned, we now use Theorem \ref{maindl} and its special form, i.e.,  Theorem \ref{coro000}, to show  all partial theta function identities (\ref{2}), (\ref{1}), and (\ref{b4}) in a unified way.

{\bf Proofs of Lemmas \ref{yl2}, \ref{yl2-12}  and \ref{yl3}.}\,\, To establish (\ref{1}), we need only to take $m=0$ in (\ref{mainid}). In such a case, we obtain at once $V_{0,n}(a,b)=1$ and
 \[U_0(b)=f(b,b)=\theta(q,b).\]
 This reduces (\ref{mainid}) to (\ref{1}).

  To achieve (\ref{b4}), recall that
     \begin{eqnarray}
     L(aq,bq)=\poq{q,aq,bq}{\infty}\sum_{n=0}^\infty\frac{\poq{ab}{2n}}{\poq{q,aq,bq,ab}{n}}q^{n}.
   \end{eqnarray}
 Meanwhile,  it is evident but important that the LHS  of (\ref{000}) is independent of $b$ while $L(a,b)$ on the RHS of (\ref{000})  is symmetric with respect to $a$ and $b$, i.e.,
   \[L(a,b)=L(b,a).\]
   Thus we  immediately obtain
     \begin{eqnarray}
 \theta(q,b)=L(a,b)+aL(aq,bq).\label{111}
   \end{eqnarray}
Subtracting (\ref{111}) from (\ref{000}) gives rise to
   \begin{eqnarray*}
    \theta(q,a)-\theta(q,b)&=&(b-a)L(aq,bq)\\
    &=&(b-a)\poq{q,aq,bq}{\infty}\sum_{n=0}^\infty
 \frac{\poq{ab}{2n}}{\poq{q,aq,bq,ab}{n}}q^n.
   \end{eqnarray*}
  It is in agreement with Schilling--Warnaar's identity (\ref{b4}).

    Finally,  to show Warnaar's identity (\ref{2}),  we need only  to  compute
    \begin{eqnarray*}
  \theta(q,a)+\theta(q,b)-1&=&\theta(q,a)-b\sum_{n=0}^\infty\tau(n)(bq)^n\\
  &=&L(a,b)+bL(aq,bq)-b\theta(q,bq).
     \end{eqnarray*}
Referring to  (\ref{111}),  we  have
\begin{eqnarray}
\theta(q,bq)=L(a,bq)+aL(aq,bq^2).\label{222}
\end{eqnarray}
Observe that $\theta(q,bq)$ is independent of $a$, which allows us to  replace $a$ with $aq$ in (\ref{222}). The result is \begin{eqnarray*}
\theta(q,bq)=L(aq,bq)+aqL(aq^2,bq^2).
\end{eqnarray*}
Therefore we are led to
 \begin{eqnarray*}
 \theta(q,a)+\theta(q,b)-1=L(a,b)-abqL(aq^2,bq^2).
  \end{eqnarray*}
 It remains to show
    \begin{eqnarray*}
  L(a,b)-abqL(aq^2,bq^2)=\poq{q,a,b}{\infty}\sum_{n=0}^\infty\frac{\poq{ab/q}{2n}}{\poq{q,a,b,ab}{n}}q^{n},
  \end{eqnarray*}
  or equivalently,
   \begin{eqnarray}
  \frac{L(a,b)}{\poq{q,a,b}{\infty}}-abq\frac{L(aq^2,bq^2)}{\poq{q,a,b}{\infty}}=
  \sum_{n=0}^\infty\frac{\poq{ab/q}{2n}}{\poq{q,a,b,ab}{n}}q^{n}.\label{mamade}
  \end{eqnarray}
  The validity of (\ref{mamade}) is ensured   by the following contiguous relation
  \begin{eqnarray*}
  &&t(a,b;n+2)=\frac{abq}{(1-a)(1-aq)(1-b)(1-bq)}t(aq^2,bq^2;n)\\
  &+&\frac{q(q-ab)(1-abq^{2n+2})}{(1-a)(1-aq)(1-b)(1-bq)(1-q^{n+1})(1-q^{n+2})}t(aq^2,bq^2;n),
  \end{eqnarray*}
  where $t(a,b;n)$ denotes the summand
  \[\frac{\poq{ab/q^2}{2n}q^{n}}{\poq{q,a,b,ab/q^2}{n}}\]
  such that
\begin{eqnarray*}
\sum_{n=0}^\infty t(a,b;n)=\frac{L(a,b)}{\poq{q,a,b}{\infty}}.
\end{eqnarray*}
Summing up, Lemmas \ref{yl2}, \ref{yl2-12}, and \ref{yl3} have been proved.
\qed

From the above derivation it is very easy to deduce that
\begin{tl} Let $\theta(q,x)$ be  given by (\ref{thetafun}). Then for any  $a,b$, it holds
 \begin{align}
  &&\frac{a}{a-b}\theta(q,a)-\frac{b}{a-b}\theta(q,b)=\poq{q,a,b}{\infty}\sum_{n= 0}^\infty \frac{\poq{ab/q^2}{2n}q^n}{\poq{q,a,b,ab/q^2}{n}}.\label{b5}\end{align}
    \end{tl}
Furthermore, with the help of  Theorem \ref{maindl-tl2}, we  are able to  show in a similar manner  that
  \begin{tl} Let $\theta(q,x)$ and $P(a,b)$ be  given by (\ref{thetafun}) and  (\ref{important-m2}),  respectively. Then
 \begin{align}
 \frac{a^2\left(b+q\right)}{
   a-b}\theta(q,a)-\frac{b^2\left(a+q
   \right)}{a-b}\theta(q,b)\nonumber\\
   =q (a+b)P(a,b)+a b(q+1)P(aq,bq).\label{b6}\end{align}
   In particular,
   \begin{align}
\left(q+a\right)\theta(q,-a)-\left(q-a\right)\theta(q,a)=2a(1+q)P(aq,-aq).\label{b7}
\end{align}
    \end{tl}
It is worth noting that  (\ref{b7}) covers a well--known result for Ramanujan's theta function $\psi(q)=\theta(q,-q)$.
\begin{lz}
     \begin{align}
\sum_{n=0}^\infty q^{n(n+1)/2}=(q^4;q^4)_{\infty}\poqq{-q}{\infty}.
\end{align}
    \end{lz}
    \pf It suffices to set $a=q$ in (\ref{b7}) and then to apply Gauss's ${}_2\phi_1$ sum (\cite[(II.8)]{10}) to $P(q^2,-q^2)$.
        \qed
\section{Further applications of $L(a,b)$}
As Eq.(2.6) of \cite{warnaar} exhibits, the LHS of (\ref{2}) is  expressed as power series in $a$ and $b$. Such expression is key to Warnaar's proof of  (\ref{2}). A quick glace at (\ref{important}) shows that $L(a,b)$ can also be expanded in a power series in $a$ and $b$. Assume now that
\begin{eqnarray}
L(a,b)&=&\sum_{i,j\geq 0}\lambda_{i,j}a^ib^j\label{symmetric}\\
&=&\poq{q,a,b}{\infty}\sum_{n=0}^\infty\frac{\poq{ab/q^2}{2n}}{\poq{q,a,b,ab/q^2}{n}}q^{n}\nonumber.
   \end{eqnarray}
Here  a  question arises naturally as follows: is there any closed expression for $\lambda_{i,j}$?  There is.
The expression  can be given by the following
\begin{dl}
Let $\lambda_{i,j}$ be defined by (\ref{symmetric}).
Then  for any nonnegative integers $i,j$, we have
\begin{eqnarray}
\lambda_{i,j}=\tau(i+j).
\end{eqnarray}
\end{dl}
\pf
From Schilling--Warnaar's identity (\ref{b4}), it follows that
\begin{eqnarray*}
L(aq,bq)&=&\poq{q,aq,bq}{\infty}\sum_{n=0}^\infty\frac{\poq{ab}{2n}}{\poq{q,aq,bq,ab}{n}}q^{n}\nonumber\\
   &=&-\sum_{n=1}^\infty\tau(n)\frac{a^n-b^n}{a-b}.
\end{eqnarray*}
 On making the substitutions  $(a,b)\to (a/q,b/q)$, we get
\begin{eqnarray*}
L(a,b)=-\sum_{n=1}^\infty\tau(n)q^{1-n}\frac{a^n-b^n}{a-b}
=\sum_{n=0}^\infty\tau(n)\frac{a^{n+1}-b^{n+1}}{a-b}.
\end{eqnarray*}
Note that
\[\frac{a^{n+1}-b^{n+1}}{a-b}=\sum_{i+j=n}a^ib^j.\]
Hence we obtain
\begin{eqnarray*}
L(a,b)=\sum_{i,j=0}^\infty\tau(i+j)a^ib^j.
\end{eqnarray*}
From here we extract  $\lambda_{i,j}=\tau(i+j),$ proving what we wanted.
\qed

As a good application of this expression of $L(a,b)$, we proceed to extend Warnaar's identity (\ref{2}) to
\begin{dl} For any integers $r,s\geq 0$ and complex numbers $a,b$, we have
\begin{eqnarray}
&&\sum_{i=0}^{r-1}\tau(i)a^i\theta(q,bq^i)+\sum_{i=0}^{s-1}\tau(i)b^i\theta(q,aq^i)\label{genralization}\\
&-&\sum_{i=0}^{r-1}\sum_{j=0}^{s-1}\tau(i+j)a^ib^j=L(a,b)-a^rb^s\tau(r+s)L(aq^{r+s},bq^{r+s}).\nonumber
\end{eqnarray}
\end{dl}
\pf It needs only to evaluate in a straightforward  way that
\begin{eqnarray*}
\mbox{RHS of (\ref{genralization}) } &=&\sum_{i,j=0}^\infty\tau(i+j)a^ib^j-a^rb^s\tau(r+s)\sum_{i,j=0}^\infty\tau(i+j)a^ib^jq^{(i+j)(r+s)}\\
&=&\sum_{i,j=0}^\infty\tau(i+j)a^ib^j-\tau(r+s)\sum_{i,j=0}^\infty\tau(i+j)a^{i+r}b^{j+s}q^{(i+j)(r+s)}\\
&=&\sum_{i=0}^{r-1}\sum_{j=0}^{\infty}\tau(i+j)a^ib^j+\sum_{i=r}^{\infty}\sum_{j=0}^{s-1}\tau(i+j)a^ib^j\\
&+&\sum_{i=r,j=s}^\infty\bigg\{\tau(i+j)-\tau(i+j-r-s)\tau(r+s)q^{(i+j-r-s)(r+s)}\bigg\}a^{i}b^{j}.
\end{eqnarray*}
As a key step,  it is easy to check that for any integers $i,j,r,s$,
\[\tau(i+j)=\tau(i+j-r-s)\tau(r+s)q^{(i+j-r-s)(r+s)}.\]
This reduces the preceding sum to
\begin{eqnarray*}
\mbox{RHS of (\ref{genralization}) } =\sum_{i=0}^{r-1}\sum_{j=0}^{\infty}\tau(i+j)a^ib^j+\sum_{i=r}^{\infty}\sum_{j=0}^{s-1}\tau(i+j)a^ib^j.
\end{eqnarray*}
This completes the proof of the theorem.
\qed

We remark that  the case $r=0,s=1$ of (\ref{genralization}) reduces to identity (\ref{000}) and the case $r=s=1$ turns out to be  (\ref{2}). Note that $\sum_{i=0}^{-1}\bullet=0.$

At the end of this paper, we employ  Warnaar's argument that used in Lemma \ref{warnaar-bailey}  to show Theorem \ref{coro111}.

{\bf Proof of Theorem \ref{coro111}.} It suffices to set $(a,b)\to (aq^{k+1},bq^{k+1}$)    simultaneously  in (\ref{important}). So we have
 \begin{eqnarray*}
 L(aq^{k+1},bq^{k+1})&=&\poq{q,aq^{k+1},bq^{k+1}}{\infty}\sum_{n=0}^\infty
 \frac{\poq{abq^{2k}}{2n}}
 {\poq{q,aq^{k+1},bq^{k+1},abq^{2k}}{n}}q^{n}\\
&\stackrel{n\to n-k}{===}&\poq{q,aq,bq}{\infty}\sum_{n=k}^\infty\frac{q^{-k}}
 {\poq{q}{n-k}\poq{ab}{n+k}}\frac{\poq{ab}{2n}}
 {\poq{aq,bq}{n}}q^{n}.
   \end{eqnarray*}
 Multiply both sides by $q^k\alpha_k\big(ab/q\big)$ and sum  the resulting identity  over all nonnegative integers $k$. Then we obtain
 \begin{eqnarray*}
 \sum_{k=0}^\infty L(aq^{k+1},bq^{k+1})q^k\alpha_k\big(ab/q\big)   =\poq{q,aq,bq}{\infty}\sum_{n=0}^\infty
  \frac{\poq{ab}{2n}}
 {\poq{aq,bq}{n}}q^{n}\beta_n\big(ab/q\big),
   \end{eqnarray*}
   where, referring back to (\ref{defnew}),  the Bailey pair relative to $ab/q$
   \[\beta_n\big(ab/q\big)=\sum_{k=0}^n\frac{\alpha_k\big(ab/q\big)}
 {\poq{q}{n-k}\poq{ab}{n+k}}.\] The theorem is thus proved.
\qed

We end our paper with a few concrete examples of Theorem \ref{coro111}.
\begin{lz}
Consider the specific Bailey pair
\begin{eqnarray*}
\left\{
  \begin{array}{ll}
   & \beta_n(ab/q)=\delta_{n,0}\quad\,\,\mbox{(Kronecker delta)} \\
   & \\
   &\alpha_n(ab/q)=\displaystyle\tau(n)\frac{1-abq^{2n-1}}{1-ab/q}
   \frac{\poq{ab/q}{n}}{\poq{q}{n}}.
  \end{array}
\right.
\end{eqnarray*}
A direct substitution of this into (\ref{111-added}) gives rise to
\begin{eqnarray*}
 \sum_{n=0}^\infty L(aq^{n+1},bq^{n+1})q^n\tau(n) \frac{1-abq^{2n-1}}{1-ab/q}\frac{\poq{ab/q}{n}}{\poq{q}{n}}  =(aq,bq,q;q)_{\infty}.
   \end{eqnarray*}
  Once combining this identity with (\ref{b5}), we easily obtain
    \begin{align}
  \sum_{n=0}^\infty q^{n}\tau(n)\left\{\frac{a}{a-b}\theta(q,aq^{1+n})+\frac{b}{b-a}\theta(q,bq^{1+n})\right\}
  \frac{\poq{ab/q}{n}}{\poq{q}{n}}\frac{1-abq^{2n-1}}{1-ab/q}=(aq,bq,q;q)_{\infty}.\label{44445555}
     \end{align}
  \end{lz}
  \begin{remark}
We remark that (\ref{44445555}) is   to some degree analogous to Alladi--Berkovich's identity (\ref{44}). It is also worthy noting that there are two  specific Bailey pairs relative to $t=0$
\begin{eqnarray*}
\left\{
  \begin{array}{ll}
   & \alpha_n(0)=\displaystyle   \frac{(-1)^nx^nq^{(n^2-n)/2}}{\poq{q}{n}}\\
   & \\
   &\beta_n(0)=\displaystyle \frac{\poq{x}{n}}{\poq{q}{n}}
  \end{array}
\right.
\quad\mbox{and}\qquad
\left\{
  \begin{array}{ll}
   & \alpha_n(0)=\displaystyle   \frac{x^nq^{n^2}}{\poq{q,xq}{n}}\\
   & \\
   &\beta_n(0)=\displaystyle \frac{1}{\poq{q,xq}{n}}.
  \end{array}
\right.
\end{eqnarray*}
Note that  for $b=0$, $L(a,0)=\theta(q,a)$. In such cases,  (\ref{111-added}) reduces respectively to
\begin{eqnarray}
 \sum_{n=0}^\infty \theta(q,aq^{n+1})\frac{(-1)^nx^nq^{(n^2+n)/2}}{\poq{q}{n}}&=&\poq{q,aq}{\infty}\sum_{n=0}^\infty
  \frac{\poq{x}{n}}
 {\poq{q,aq}{n}}q^{n},\\\label{111-added-added-added}
   \sum_{n=0}^\infty \theta(q,aq^{n+1})\frac{x^nq^{n^2+n}}{\poq{q,xq}{n}}&=&\poq{q,aq}{\infty}\sum_{n=0}^\infty
  \frac{q^{n}}
 {\poq{q,aq,xq}{n}}.\label{111-added-added}
   \end{eqnarray}
   \end{remark}
\section*{Acknowledgements}
 This  work was supported by the National Natural Science Foundation of China [Grant No.  11471237].

 \bibliographystyle{amsplain}

\end{document}